\documentclass[10pt]{article}
\usepackage{amssymb,epsfig,psfrag}
\usepackage{amsfonts, amssymb}
\usepackage{amsmath}
\usepackage{amsthm}
\usepackage{enumerate}
\usepackage{amsmath, amssymb, amsthm,algorithm}
\usepackage{tocloft}

\newtheorem{thm}{Theorem}[subsection]

\newtheorem{lemma}[thm]{Lemma}

\theoremstyle{remark}
\newtheorem*{rem}{Remark}
\newtheorem*{ex}{Example}

\begin{document}

\title{Periodicity of Multidimensional Continued Fractions\\} 
\author{Eun Hye Lee}
\date{}
\maketitle

\begin{abstract}
It is known that the continued fraction expansion of a real number is periodic if and only if the number is a quadratic irrational. In an attempt to generalize this phenomenon to other settings, Jun-Ichi Tamura and Shin-Ichi Yasutomi have developed a new algorithm for multidimensional continued fractions (Algebraic Jacobi-Perron algorithm) that involves cubic irrationals, and proved periodicity in some cubic number fields, such as $\mathbb{Q}(\sqrt[3]{m^3+1})$ where $m\in\mathbb{Z}$, and $\mathbb{Q}(\delta_m)$ where $\delta_m$ is a root of $x^3-mx+1=0,\,\,m\in\mathbb{Z},\,\, m\geq3$ with the algorithm.

In this paper, we study some other types of number fields that give rise to periodic continued fractions using the Algebraic Jacobi-Perron algorithm obtaining results for $\mathbb{Q}(\sqrt[l]{m^l+1})$ for any positive integer $l$.
Furthermore, we find that some families of cubic equations, such as  $x^3+3ax^2+bx+ab-2a^3+1=0,\,b\leq3a^2-3,\,a,b\in\mathbb{Z}$, have roots that have periodic multidimensional continued fractions.

\begin{itemize}
\item \textbf{Key Words:} Jacobi-Perron algorithm, Multidimensional continued fraction, Periodicity
\end{itemize}
{\bf Aknowledgement} This paper is due to my master's degree thesis. I would like to express my gratitude for Dr. Yoonjin Lee, my mentor and Master advisor, for her continued guidance and encouragement through my degree. I would also like to thank Dr. Ramin Takloo-Bighash, my Ph.D advisor, for his support in getting this paper to see the light of the day.
\end{abstract} 
\pagenumbering{arabic}
\pagestyle{myheadings}

\section{Introduction}

One of the most well-known properties of continued fractions is the result that the simple continued fraction expansion of a real number is periodic if and only if that real number is a quadratic irrational. In order to obtain other periodicity results of this nature mathematicians have had to generalize the types of continued fractions they have considered. 
The types of continued fractions we consider in this paper are {\em multidimensional} in the sense we will make precise later.  Our multidimensional continued fractions are produced using the Jacobi-Perron type algorithms which have previously been studied by Jun-Ichi Tamura and Shin-Ichi Yasutomi. In \cite{MR1183919}, Tamura introduces the Jacobi-Perron Algorithm and then in \cite{MR1339134}, \cite{MR2521286}, and \cite{MR2731821}, Tamura and Yasutomi discuss modified versions of Jacobi-Perron Algorithm, such as the Jacobi-Perron-Parusnikov Algorithm and Algebraic Jacobi-Perron Algorithm, and their connections with higher dimensional continued fractions. In particular, in \cite{MR2521286}, Tamura proved that certain types of elements in number fields can be expanded as periodic continued fractions. In addition, Adam \cite{MR3612875} and Voutier \cite{MR3515830} have investigated the periodicity of Jacobi-Perron algorithms.

	In this paper, we study some other types of number fields whose elements  can be expressed as periodic continued fractions using the Algebraic Jacobi-Perron algorithm. Tamura proved such a result for $\mathbb{Q}(\sqrt[3]{m^3+1})$. Here we generalize the base field to be $\mathbb{Q}(\sqrt[l]{m^l+1})$ for any $l\ge2$ using similar method Tamura used but for larger degree through the observations of how the terms shift through the algorithm. Furthermore, we find that some families of cubic equations, such as  $x^3+3ax^2+bx+ab-2a^3+1=0,\,b\leq3a^2-3,\,a,b\in\mathbb{Z}$, have roots that become periodic when expanded as continued fractions using the same algebraic Jacobi-Perron algorithm. This uses the result of Tamura, but using the general facts of norms of cubic number fields, we get more general result.
\bigskip

Our first theorem is the following: 
\begin{thm} For natural numbers $m, l>1$, the elements 
$$
\sqrt[l]{(m^l+1)^{k}}-m^{k}, \quad 0 < k < l
$$
of $\mathbb{Q}(\sqrt[l]{m^l+1})$ have periodic continued fractions. 
\end{thm}

\begin{ex} Consider the case when $l=4$. The vector 
$\begin{bmatrix}\sqrt[4]{m^4+1}-m\\ \sqrt[4]{(m^4+1)^2} -m^2\\ \sqrt[4]{(m^4+1)^3}-m^3\end{bmatrix} $ is periodic with period 3 for any $m>1$.
\begin{center}
  \begin{tabular}{ c | c  c  c  c  c c}
   
    $n$ & $1$ & $2$ & $3$ & $4$ & $5$ & $n_{(\geq6)}$ \\ \hline
    $a_n$ & $4m^3$ & $6m^2$ & $4m$ & $4m^3$ & $6m^2$ & $a_{n-3}$\\ 
    $b_n$ & $2m$ & $4m^3$ & $6m^2$ & $4m$ & $4m^3$ & $b_{n-3}$\\
    $c_n$ & $3m^2$ & $3m$ & $4m^3$ & $6m^2$ & $4m$ & $c_{n-3}$
  \end{tabular}
\end{center}
\end{ex}
\bigskip

The following is our second theorem: 
\begin{thm}\label{thm2} Suppose $a, b$ are integers, and let $\gamma$ be a root of $x^3+3ax^2+bx+ab-2a^3+1=0$. Furthermore, assume $-a<\gamma<-a+1$. Then the real numbers 
$$\gamma-\lfloor\gamma\rfloor, \quad (\gamma-\lfloor\gamma\rfloor)^2$$ 
have periodic continued fraction expansions. 
\end{thm}


\begin{ex} Consider the function $f(x)=x^3+3x^2-1$, i.e. $a=1$ and $b=0$ in \ref{thm2}. Then there is a solution $\gamma$ in between $-a=-1$ and $-a+1=0$. In fact, $\gamma\approx -0.65270$. Now, using the transformation, we know that $\delta=\gamma+1=\gamma-\lfloor\gamma\rfloor$ is a solution to $g(x)=x^3-3x+1$. The vector $\begin{bmatrix}\gamma-\lfloor\gamma\rfloor\\(\gamma-\lfloor\gamma\rfloor)^2 \end{bmatrix}$ is periodic with period 4.
\begin{center}
  \begin{tabular}{ c | c  c  c  c  c  c  c  c}
   
    $n$ & $1$ & $2$ & $3$ & $4$ & $5$ & $6$ & $7$ & $n_{(\geq8)}$ \\ \hline
    $a_n$ & $2$ & $1$ & $0$ & $0$ & $1$ & $-1$ & $1$ & $a_{n-4}$\\ 
    $b_n$ & $0$ & $0$ & $2$ & $-1$ & $1$ & $0$ & $1$ & $b_{n-4}$
  \end{tabular}
\end{center}

\end{ex}
\bigskip
There are numerous applications to normal continued fractions, such as cryptography and approximation. However, applications for these multidimensional continued fractions are not yet found. 

\bigskip

This is how this paper is organized. In \S \ref{algebraic} we review the Algebraic Jacobi-Perron algorithm and discuss some previous works. \S \ref{our results} contains  our main results.

\subsection*{Notation}

The following notations are used through the paper.\\
$K=\mathbb{C}(z^{-1})$ : the field of formal Laurent series with complex coefficients.\\
$||\varphi||:=e^{-k}$ for $\displaystyle0\neq \varphi = \sum^\infty_{m=k} c_mz^{-m}\in K$, $c_k\neq0$, $k\in \mathbb{Z}$.\\
$||0||:=0$.\\
${\rm{d}}_K(\varphi, \psi):=||\varphi-\psi||$ $(\varphi, \psi \in K)$ : the distance function.\\
$[\varphi]:=$ the polynomial part of $\varphi\in K$.\\
($i.e.$ $ [\varphi]:=0$ for $||\varphi||<1$ and $[\varphi]:=\sum^0_{m=k} c_mz^{-m}$ for $||\varphi||\geq1$ with $\varphi=\sum^\infty_{m=k} c_mz^{-m}$)\\
$<\varphi>:=\varphi-[\varphi]$ : fractional part of $\varphi\in K$.\\
For $\boldsymbol{\varphi}\in K^l$, \\
$\boldsymbol{\varphi}:=\left(\begin{matrix} \varphi_1\\ \varphi_2\\ \vdots \\\varphi_l \end{matrix}\right)$,\,\,\,
$[\boldsymbol{\varphi}]:=\left(\begin{matrix} [\varphi_1]\\ [\varphi_2]\\ \vdots \\ [\varphi_l ] \end{matrix}\right)$, $<\boldsymbol{\varphi}>:=\left(\begin{matrix} <\varphi_1>\\ <\varphi_2>\\ \vdots \\<\varphi_l> \end{matrix}\right)$.\\

$T\left(\begin{matrix} \varphi_1\\ \varphi_2\\ \vdots \\\varphi_l \end{matrix}\right):=\left(\begin{matrix} {1\over\varphi_l}\\ {\varphi_1\over\varphi_l}\\ \vdots \\{\varphi_{l-1}\over\varphi_l} \end{matrix}\right):={\cfrac{1}{\left(\begin{matrix} \varphi_1\\ \varphi_2\\ \vdots \\\varphi_l \end{matrix}\right)}}$ : transformation.\\

$\mathcal{P}_K:=$set of elements in $K$ which have periodic continued fractions.
\bigskip

\bigskip


\section{Algebraic Jacobi-Perron Algorithm}\label{algebraic}

In this section we recall several constructions of multidimensional continued fraction expansions of real numbers. 

\subsection{Jacobi-Perron-Parusnikov Algorithm}
\paragraph{}
Put \begin{equation*}\boldsymbol{\varphi_0}=\left(\begin{matrix} \varphi_{1_0}\\ \varphi_{2_0}\\ \vdots \\\varphi_{l_0} \end{matrix}\right):=\left(\begin{matrix} <\varphi_1>\\ <\varphi_2>\\ \vdots \\<\varphi_l> \end{matrix}\right).\end{equation*}
Identify $$\boldsymbol{\varphi_n}\equiv \left(\begin{matrix} \varphi_{1_n}\\ \varphi_{2_n}\\ \vdots \\\varphi_{l_n} \end{matrix}\right):=S^n(\boldsymbol{\varphi_0}) \mbox{ where }S(\boldsymbol{\varphi_0}):=<T^{-1}(\boldsymbol{\varphi_0})>.$$
 ($S^n$ : the $n$-fold iteration of $S$).\\
 Then we define $$\bf{a_n}\equiv \left(\begin{matrix} a_n\\ b_n\\ \vdots \\l_n \end{matrix}\right)\in (\mathbb{C}[z])^2$$ by
 \begin{equation}\label{eq:def1} {\bf{a_0}}:=[\boldsymbol{\varphi}], \;{\bf{a_n}}:=[T^{-1}(S^{n-1}(\boldsymbol{\varphi_0}))] \;(n\geq1).\end{equation}
 Then the vector $\pi_n\in(\mathbb{C}(z))^2$ defined by
 \begin{equation}\label{eq:def2} \pi_n=\pi_n(\boldsymbol{\varphi})=\bf{a_0} + \cfrac{1}{\bf{a_1} + \cfrac{1}{\ddots \, + \cfrac{1}{\bf{a_n}}}}\end{equation}
 converges to $\boldsymbol{\varphi}$ component wise as $n\rightarrow\infty$ (Parysnikov's Thm) where the convergence is with respect to the metric d$_K$. Thus, we can write
\begin{equation}\label{eq:def3} \pi=\bf{a_0} + \cfrac{1}{\bf{a_1} + \cfrac{1}{\ddots \, + \cfrac{1}{\bf{a_n} + \ddots}}} \end{equation}
(\ref{eq:def2}) and (\ref{eq:def3}) will be also written as 
\begin{equation*} \pi_n\left(\begin{matrix} \varphi_1\\ \varphi_2\\ \vdots \\\varphi_l \end{matrix}\right)=\left[\begin{matrix} a_0;&a_1,&a_2,&\dots,&a_n\\b_0;&b_1,&b_2,&\dots,&b_n\\\vdots\\ l_0;&l_1,&l_2,&\dots,&l_n\end{matrix}\right],\end{equation*}
\begin{equation*} \left(\begin{matrix} \varphi_1\\ \varphi_2\\ \vdots \\\varphi_l \end{matrix}\right)=\left[\begin{matrix} a_0;&a_1,&a_2,&\dots,&a_n,&\dots\\b_0;&b_1,&b_2,&\dots,&b_n,&\dots\\\vdots\\ l_0;&l_1,&l_2,&\dots,&l_n,&\dots\end{matrix}\right].\end{equation*}

The right-hand side of (\ref{eq:def3}) is called "the expression of $\boldsymbol{\varphi}\in K^l$ in the {\bf{\textit{Jacobi-Perron-Parusnikov Algorithm (JPPA)}}}".\newline

We can consider the expression (\ref{eq:def3}) for any given sequence of vectors $\bf{a_1,\,a_2,\,\dots}$ in $K^l$ provided that its $n^{th}$ convergent is well-defined and converges to some element of $K^l$. We say that the JPPA is admissible if it is derived from (\ref{eq:def1}).\bigskip

{\ex{If we take the field $\mathbb{R}$ of real numbers for $K$ and define $[\varphi]$ to be the integral part of $\varphi\in\mathbb{R}$, we have algorithm which is the simplest JP algorithm.}}

Let $$p(x)=x^3-kx^2-lx-1$$ be a polynomial satisfying $$k\geq l\geq0, k+l\geq2, (k,\,l\in\mathbb{Z}).$$ Let $\lambda=1/\alpha$ and $\kappa=l/\alpha+a/\alpha^2$ where $\alpha>1$ is the root of $p(x)$. Then the following expression is admissible in the JP algorithm.
\begin{equation}\label{eq:ex} \left(\begin{matrix} \lambda\\\kappa\end{matrix}\right)=\left[\begin{matrix} 0;&l,&l,&l,&\dots\\0;&k,&k,&k,&\dots\end{matrix}\right] \end{equation}
\begin{proof} We have $T^{-1}\left(\begin{matrix} \lambda\\\kappa\end{matrix}\right)=\left(\begin{matrix} l+1/\alpha\\\alpha\end{matrix}\right)=\left(\begin{matrix}l\\k\end{matrix}\right)+\left(\begin{matrix}\lambda\\\kappa\end{matrix}\right)$ with $0<\lambda<1$ and $0<\kappa<1$ which leads to the admissible expression (\ref{eq:ex}).

\end{proof}
\bigskip

\subsection{Algebraic Jacobi-Perron Algorithm (degree 3)}

Let $K$ be a number field of degree 3.
We denote $X_K$ as the set defined by
$$
X_K:=\{(\alpha,\beta)\in K^2 | 1,\alpha,\beta \mbox{ are linearly independent over }\mathbb{Q}\}\cap I^2$$ where $ I=[0,1).$
We define the transformation $T_K$ on $X_K$ by

\begin{equation*}
T_K(\alpha,\beta):=\left\{\begin{array}{rl}
\left(\cfrac{1}{\alpha}-\left\lfloor\cfrac{1}{\alpha}\right\rfloor,\cfrac{\beta}{\alpha}-\left\lfloor\cfrac{\beta}{\alpha}\right\rfloor\right) & \mbox{ if }\cfrac{\alpha}{\sqrt{|N(\alpha)|}}>\cfrac{\beta}{\sqrt{|N)\beta)|}}, \\
\left(\cfrac{\alpha}{\beta}-\left\lfloor\cfrac{\alpha}{\beta}\right\rfloor, \cfrac{1}{\beta}-\left\lfloor\cfrac{1}{\beta}\right\rfloor\right) & \mbox{ if }\cfrac{\alpha}{\sqrt{|N(\alpha)|}}<\cfrac{\beta}{\sqrt{|N)\beta)|}}\end{array} \right.
\end{equation*}
for $(\alpha,\beta)\in X_K$, where $\lfloor x \rfloor$ is the floor function of $x$ and $N(x)$ is the norm of $x\in K$ over $\mathbb{Q}$. Also, define integral valued functions $a(\alpha,\beta)$ and $b(\alpha,\beta)$,$$a(\alpha,\beta)=\begin{cases}
	\left\lfloor \cfrac{1}{\alpha}\right\rfloor, &\mbox{ if }\cfrac{\alpha}{\sqrt{|N(\alpha)|}}>\cfrac{\beta}{\sqrt{|N(\beta)|}}\\
	\left\lfloor \cfrac{\alpha}{\beta}\right\rfloor, &\mbox{ if }\cfrac{\alpha}{\sqrt{|N(\alpha)|}}<\cfrac{\beta}{\sqrt{|N(\beta)|}}\end{cases}$$
	
	$$b(\alpha,\beta)=\begin{cases}
	\left\lfloor \cfrac{\beta}{\alpha}\right\rfloor, &\mbox{ if }\cfrac{\alpha}{\sqrt{|N(\alpha)|}}>\cfrac{\beta}{\sqrt{|N(\beta)|}}\\
	\left\lfloor \cfrac{1}{\beta}\right\rfloor, &\mbox{ if }\cfrac{\alpha}{\sqrt{|N(\alpha)|}}<\cfrac{\beta}{\sqrt{|N(\beta)|}}\end{cases}$$

Finally, the continued fraction form of $(\alpha,\beta)$ becomes
	$$(\alpha,\beta)=\left[\begin{matrix}
	a(\alpha,\beta);&a(\alpha_1,\beta_1),&a(\alpha_2,\beta_2)&\dots\\
	b(\alpha,\beta);&b(\alpha_1,\beta_1),&b(\alpha_2,\beta_2)&\dots
	\end{matrix}\right].$$

\bigskip

\subsection{Previous Works over Some Cubic Number Fields}

\bigskip

The following are theorems from \cite{MR2521286}, which inspired this paper. They show the periodicity of elements of certain cubic number fields.

\bigskip

{\thm{[Theorem 2.4. \cite{MR2521286}]}\label{pure} Let $K=\mathbb{Q}(\sqrt[3]{m^3+1})$ with $m\in\mathbb{Z}_{>0}$. Moreover, let$$(\alpha,\beta)=(\sqrt[3]{m^3+1}-m,\sqrt[3]{(m^3+1)^2}-m^2).$$ Then, $(\alpha,\beta)\in \mathcal{P}_K$ and the length of the period is 2.}

\begin{itemize}
\item{$m\geq2$}
\begin{center}
  \begin{tabular}{ c | c  c  c  c}
   
    $n$ & $1$ & $2$ & $3$ & $n_{(\geq4)}$ \\ \hline
    $a_n$ & $2m^2$ & $3m$ & $3m^2$ & $a_{n-2}$\\ 
    $b_n$ & $2m$ & $3m^2$ & $3m$ & $b_{n-2}$ 
  \end{tabular}
\end{center}

\item{$m=1$}
\begin{center}
\begin{tabular}{c | c  c  c}
$n$ & $1$ & $2$ & $n_{(\geq3)}$ \\ \hline
$a_n$ & 0 & 2 & $a_{n-2}$\\
$b_n$ & 1 & 1 & $b_{n-2}$
\end{tabular}
\end{center}
\end{itemize}
\bigskip

{\thm{[Theorem 2.5. \cite{MR2521286}}]\label{sec:cubic} Let $\delta_m$ be the root of $x^3-mx+1=0$ $(m\in \mathbb{Z}, m\geq3)$ determined by $0<\delta_m<1$. Then, $K=\mathbb{Q}(\delta_m)$ is a cubic number field and $(\delta_m, \delta_m^2)\in \mathcal{P}_K$ with the period 4.}

\begin{center}
  \begin{tabular}{ c | c  c   c  c  c  c  c  c}
   
    $n$ & $1$ & $2$ & $3$ & $4$ & $5$ & $6$ & $7$ & $n_{(\geq8)}$ \\ \hline
    $a_n$ & $m-1$ & $1$ & $0$ & $0$ & $m-2$ & $1$ & $1$ & $a_{n-4}$\\ 
    $b_n$ & $0$ & $0$ & $m-1$ & $1$  & $1$ & $0$ & $m-2$ & $b_{n-4}$
  \end{tabular}
\end{center}

\bigskip
\section{Multidimensional Continued Fractions} \label{our results}

\vspace{5mm}
\hspace{7mm}
In this chapter, we see how to make multidimensional continued fractions using Algebraic Jacobi-Perron algorithm. Then we discuss about the periodicity of the continued fraction obtained by it.

\subsection{Algebraic Jacobi-Perron Algorithm (higher degree)}
\paragraph{}

Let $K$ be a number field of degree $l$ and let $$X_K:=\{(\alpha_1,\dots,\alpha_{l-1})\in K^{n-1}|1, \alpha_i\mbox{ are linearly independent over }\mathbb{Q}\}\cap I^{n-1}$$ where $I=[0,1)$.

We define the integer-valued functions $a_1$, $a_2$, \dots, $a_{l-1}$ on $X_K$ as follows:
$$
a_i(\alpha_1,\alpha_2,\dots,\alpha_{l-1}) =\begin{cases}
 \left\lfloor \cfrac{1}{\alpha_i} \right\rfloor &\mbox{ if $ \cfrac{\alpha_i}{\sqrt[l-1]{|N(\alpha_i)|}} > \cfrac{\alpha_j}{\sqrt[l-1]{|N(\alpha_j)|}}$}\\
&\hfill\mbox{ for $ 1\leq j\leq l-1, i\neq j$}\\
 \left\lfloor \cfrac{\alpha_i}{\alpha_j}\right\rfloor &\mbox{ for the index of the biggest $\cfrac{\alpha_j}{\sqrt[l-1]{|N(\alpha_j)|}}$.}
 \end{cases}   
$$
for $(\alpha_1,\alpha_2,\dots,\alpha_{l-1})\in X_K$.

We define the transformation $T_K$ on $X_K$ by
\begin{equation*} T_K(\alpha_1,\dots,\alpha_{l-1}):=(b_1(\alpha_1,\dots,\alpha_{l-1}), \dots, b_{l-1}(\alpha_1,\dots, \alpha_{l-1})) \end{equation*}
where
$$ b_i(\alpha_1,\alpha_2,\dots,\alpha_{l-1}) =  \begin{cases}
\cfrac{1}{\alpha_i}- \left\lfloor \cfrac{1}{\alpha_i} \right\rfloor &\mbox{ if $\cfrac{\alpha_i}{\sqrt[l-1]{|N(\alpha_i)|}} > \cfrac{\alpha_j}{\sqrt[l-1]{|N(\alpha_j)|}} $}\\&\hfill \mbox{  for $1\leq j\leq l-1$, $i\neq j$} \\
 \cfrac{\alpha_i}{\alpha_j}-\left\lfloor \cfrac{\alpha_i}{\alpha_j}\right\rfloor &\mbox{ for the index of the biggest $\cfrac{\alpha_j}{\sqrt[l-1]{|N(\alpha_j)|}}$\,.}
       \end{cases}$$

For $n\in \mathbb{Z}^+$, we put
\begin{multline*}
(a_{1_n}, a_{2_n},\dots,a_{{l-1}_n})  =(a_{1_n}(\alpha_1,\alpha_2, \dots,\alpha_{l-1}), a_{2_n}(\alpha_1,\alpha_2,\dots,\alpha_{l-1}),\\
\hfill\dots,a_{l-1_n}(\alpha_1,\alpha_2,\dots,\alpha_{l-1}))\,\,\,\,\\
\hfill:=(a_1(T^{n-1}_K(\alpha_1,\alpha_2,\dots,\alpha_{l-1})),a_2(T^{n-1}_K(\alpha_1,\alpha_2,\dots,\alpha_{l-1})),\,\,\,\,\\
\hfill \dots,a_{l-1}(T^{n-1}_{K}(\alpha_1,\alpha_2,\dots,\alpha_{l-1}))),
\end{multline*}
and
\begin{multline*}
S(\alpha_1,\dots, \alpha_{l-1} ) := \{(a_{1_n}(\alpha_1,\dots,\alpha_{l-1}), a_{2_n}(\alpha_1,\dots,\alpha_{l-1}),\\\hfill\dots,a_{l-1_n}(\alpha_1,\dots,\alpha_{l-1}))\}^\infty_{n=1}.
\end{multline*}
The sequence $S(\alpha_1,\dots,\alpha_{l-1})$ will be referred to as the expansion of $$(\alpha_1,\dots,\alpha_{l-1})\in X_K$$ by $T_K$; $T_K$ gives rise to a $l-1$-dimensional continued fraction expansion.
Finally, we have the continued fraction form
$$\left[\begin{matrix}\alpha_1\\\alpha_2\\\vdots\\\alpha_{n-1}\end{matrix}\right]
=\left[\begin{matrix}a_{1_1};& a_{1_2},&\dots,&a_{1_m},&\dots\\ a_{2_1},&a_{2_2},&\dots,&a_{2_m},&\dots\\\vdots&&\vdots&&\vdots\\ a_{n-1_1};&a_{n-1_2},&\dots,&a_{n-1_m},&\dots\end{matrix}\right].$$ which will be called the {\textit{\textbf{Algebraic Jacobi-Perron Alpgorithm(AJPA)}}}.
\bigskip

We show that the transformation $T_K$ is well defined.
\bigskip
{\lemma 

The transformation $T_K$ is well defined.}

\begin{proof} Let $(\alpha_1, \dots, \alpha_{l-1})\in X_K$. \newline It suffices to show that $$\cfrac{\alpha_i}{\sqrt[l-1]{|N(\alpha_i)|}} \neq \cfrac{\alpha_j}{\sqrt[l-1]{|N(\alpha_j)|}}\mbox{ for }1\leq i \neq j\leq l-1.$$ We suppose $\cfrac{\alpha_i}{\sqrt[l-1]{|N(\alpha_i)|}} = \cfrac{\alpha_j}{\sqrt[l-1]{|N(\alpha_j)|}}$. Then, we have $\alpha_i = \sqrt[l-1]{\cfrac{|N(\alpha_i)|}{|N(\alpha_j)|}}\,\alpha_j$. Since $\alpha_i$ and $\alpha_j$ are linearly independent over $\mathbb{Q}$, $\sqrt[l-1]{\cfrac{|N(\alpha_i)|}{|N(\alpha_j)|}} \notin \mathbb{Q}$. Hence $\sqrt[l-1]{\cfrac{|N(\alpha_i)|}{|N(\alpha_j)|}}$ is a $l-1$th irrational and $\sqrt[l-1]{\cfrac{|N(\alpha_i)|}{|N(\alpha_j)|}}\in K$, which is a contradiction. It is easy to see that $T_K\left(\sqrt[l-1]{\cfrac{|N(\alpha_i)|}{|N(\alpha_j)|}}\,\right)\in X_K$.

 \end{proof}
\bigskip

{\rem The transformaion $T_K$ can be applied generally to any elements of $X_K$, where $K$ is any number field. Thus, the continued fraction is obtained from any elements, whether they are periodic or not. In this thesis, we concentrate on finding periodic elements.}

\bigskip
\subsection{Periodicity of Multidimensional Continued Fractions}

Once an element is written as a continued fraction, if it is periodic, more information can be obtained. It is previously shown in 3.3.1 and 3.3.2 that
\begin{equation}\label{cubic}
(\sqrt[3]{m^3+1}-m,\sqrt[3]{(m^3+1)^2}-m^2)
\end{equation}
and
\begin{equation}\label{delta}
(\delta_m, \delta_m^2)
\end{equation}
where $\delta_m$ is a root of $x^3-mx+1=0$, $m\geq3$ and $m\in\mathbb{Z}$, with $0<\delta_m<1$ are periodic.

\bigskip

Now, in the sense of generalizaing (\ref{cubic}), we show that
$$
(\sqrt[l]{m^l+1}-m,\sqrt[l]{(m^l+1)^2}-m^2,\dots,\sqrt[l]{(m^l+1)^{l-1}}-m^{l-1})
$$
is also periodic.
\bigskip
{\thm
{Let $K=\mathbb{Q}(\sqrt[l]{m^l+1})$ with $m\in\mathbb{Z}_{\geq2}$ and let $$(\alpha_1,\alpha_2,\dots,\alpha_{l-1})=(\sqrt[l]{m^l+1}-m,\sqrt[l]{(m^l+1)^2}-m^2,\dots,\sqrt[l]{(m^l+1)^{l-1}}-m^{l-1}).$$ Then $(\alpha_1,\alpha_2,\dots,\alpha_{l-1})\in\mathcal{P}_K$ and the length of the period is $l-1$.}}

\begin{proof} Observe that for $\sqrt[l]{(m^l+1)^k}-m^k$, when $k=1$, the norm becomes minimal, $i.e.\;1$. For convenience, I shall denote $\sqrt[l]{m^l+1}$ as $x$, $i.e.$, every element of the tuple is expressed in terms of $x$ and $m$. Note that since coefficients of elements of a tuple is 0, all elements are divisible by $x-m$. Also, since all tuples have $x-m$ as one of the elements and only one of them, it always goes to denominator after taking the map $T_K$.

After taking $T_K$, the order of all the elements get decreased by 1 except for the element which is $x-m$. Taking $T_K$ to $x-m$, the element becomes $x^{l-1}+m\,x^{l-2}+\dots+m^{l-2}\,x+m^{l-1}-l\,m^{l-1}$ (which is still divisible by $x-m$).

After an iteration, we observe that only the position of the $x-m$ and elements left to $x-m$ is shifted, $i.e.$, $$T_K(a_1,\dots,a_{i-1},x-m,a_{i+1},\dots,a_{l-1})=(b_1,a_1,\dots,a_{i-1},x-m,\dots).$$ Thus, after $l-1$ iterations, all elements has been $x-m$ once, and now, the first element is $x-m$ again. So from now on, taking $T_K$ means just shifting the tuple.

That is to say, $$T_K^{l-1}(a_1,a_2,\dots,a_{l-1})=(a_1,a_2,\dots,a_{l-1})$$ where $(a_1 ,a_2, \dots, a_{l-1}) = T_K^k (\alpha_1, \alpha_2, \dots, \alpha_{l-1})$, $k\geq l-1$ and that $$(\alpha_1, \alpha_2, \dots, \alpha_{l-1})\in \mathcal{P}_K.$$\end{proof}

\bigskip

For the generalization of (\ref{delta}), first we see a lemma showing a family of cubic polynomial that can be reduced to $x^3-mx+1$ with $m\geq3$, $m\in\mathbb{Z}$.
\bigskip

{\rem{

We try to make the constant term of reduced polynomial be 1 for the convenience of computing norms. It is not easy to generalize the polynomial $x^3-mx+1$ with $m\geq3, m\in\mathbb{Z}$  in terms of its constant term because of the complexicy of computation for the norm. It is still an open problem whether or not there is more effective way to compute norms with polynomial with constant term different from 1.}}

\bigskip

{\lemma

A family of cubic polynomials $x^3+3ax^2+bx-ab-2a^3+1$ with $b\leq3a^2-3$ and $a,b\in\mathbb{Z}$ can be reduced to $x^3-(3a^2-b)x+1$.}

\begin{proof}
A cubic monic polynomial
\begin{equation}
x^3+kx^2+lx+n
\end{equation}
is reduced to
\begin{multline*}
\left(x-\cfrac{k}{3}\right)^3+k\left(x-\cfrac{k}{3}\right)^2+l\left(x-\cfrac{k}{3}\right)+n\\=x^3-\left(l-\cfrac{k^2}{3}\right)x+\left(n-\cfrac{kl}{3}+\cfrac{2k^3}{27}\right).
\end{multline*}
We need

\begin{enumerate}[(i)]
\item
$\cfrac{k^2}{3}-l\in\mathbb{Z}$, 
\item
$\cfrac{k^2}{3}-l\geq3$, and
\item
$ n-\cfrac{kl}{3}+\cfrac{2k^3}{27}=1$.
\end{enumerate}
From (i), we can conclude that $k=3a$ for some $a\in\mathbb{Z}$.
\\From (ii), $3a^2-l\geq3$, thus $l\leq3a^2-3$.
\\From (iii), $n-al+2a^3=1$, thus $n=al-2a^3+1$.

Reindexing $l$ as $b$, (7) becomes
$$
x^3+3ax^2+bx+ab-2a^3+1
$$
and it is reduced to 
$$x^3-(3a^2-b)x+1$$
as desired.

\end{proof}
\bigskip
Now, we are to show that periodic continued fraction using the root of
$$x^3+3ax^2+bx+ab-2a^3+1$$
is equivalent to see the roots of the reduced polynomial, $$x^3-(3a^2-b)x+1.$$

\bigskip

{\thm\label{sec:red}

Let $K=\mathbb{Q}(\gamma)$ where $\gamma$ is a root of $x^3+3ax^2+bx+ab-2a^3+1=0$ when $b\leq 3a^2-3$, $a,b\in\mathbb{Z}$ with $-a<\gamma<-a+1$ and also let
$$(\alpha_1,\alpha_2)=(\gamma-\lfloor\gamma\rfloor,(\gamma-\lfloor\gamma\rfloor)^2).$$
Then $(\alpha_1,\alpha_2)\in\mathcal{P}_K$ and the length of the period is 4.}

\begin{proof}

Let $\delta$ be a root of $x^3-(3a^2-b)x+1=0, b\leq3a^2-3, a,b\in\mathbb{Z}$ such that $0<\delta<1$. Then, since the transformation is linear $(i.e. x\rightarrow x+a)$, a root of $x^3+3ax^2+bx+ab-2a^3+1=0$ must be $\delta-a(=:\gamma)$. Because $a\in\mathbb{Z}$, $\delta$ and $\gamma$ have the same fractional part. That is to say, $\gamma-\lfloor\gamma\rfloor=\delta-\lfloor\delta\rfloor=\delta$. Hence,
$$(\gamma-\lfloor\gamma\rfloor,(\gamma-\lfloor\gamma\rfloor)^2)=(\delta,\delta^2).$$
We know that $(\delta,\delta^2)$ is periodic form \ref{sec:cubic}. Thus $(\gamma-\lfloor\gamma\rfloor,(\gamma-\lfloor\gamma\rfloor)^2)$ is periodic.

\end{proof}

{\rem

Every cubic polynomial can be reduced to a trinomial. \ref{sec:red} is possible since the original cubic polynomial and the reduced on have the same splitting field over $\mathbb{Q}$. Since the splitting fields are the same, they produce the same norm with the same element, so Algebraic Jacobi-Perron algorithm can be applied in the same way. However, the esssential property for the multidimensional continued fraction to be periodic is not yet determined. In this step, we just look for the ones that are periodic.}

\nocite{*}
\bibliographystyle{plain}
\bibliography{lee-thesis-copy}
\end{document}